\documentclass[12pt,leqno,fleqn,epsfig]{article}
\usepackage{amssymb, epsfig, amsmath, amsthm}
\usepackage{mathrsfs}       
     \usepackage{color}
\newcommand{\vare}{\varepsilon}
\newcommand{\R}{\mathbb{R}}

\newcommand{\const}{\operatorname*{const}}

\newcommand{\bb}{\begin{equation}}
\newcommand{\ee}{\end{equation}}
\newcommand{\bq}{\begin{eqnarray}}
\newcommand{\eq}{\end{eqnarray}}
\newcommand{\bqn}{\begin{eqnarray*}}
\newcommand{\eqn}{\end{eqnarray*}}

\newcommand{\var}{\varepsilon}
\renewcommand{\o}{\omega}

\newcommand{\intl}{\int\limits}

\newcommand{\Beweisende}{\rule{0.2cm}{0.2cm}}

\newcounter{secnum}

\newtheorem{thm}{Theorem}[section]

\theoremstyle{definition}

\newtheorem{defin}[thm]{Definition}
\newtheorem{rem}[thm]{Remark}

\title{On the regularity of  solutions to the  2D  Boussinesq  equations satisfying  Type I  conditions } 
 
\author{Dongho Chae$^*$  and J\"{o}rg Wolf $^\dagger$\\
\ \\
 Department of Mathematics\\
Chung-Ang University\\
 Seoul 156-756, Republic of Korea\\
($*$)e-mail: dchae@cau.ac.kr \\
($^\dagger$)e-mail: jwolf2603@cau.ac.kr}
\date{}
\begin{document}
\maketitle
\begin{abstract}
We prove continuation in time of the local smooth solutions  satisfying various  Type I  conditions for the 2D inviscid Boussinesq equations.
  \noindent{\bf
keywords:} Boussinesq equations, Type I condition, non blow-up

\end{abstract}

\section{Introduction}
\label{sec:-1}
\setcounter{secnum}{\value{section} \setcounter{equation}{0}
\renewcommand{\theequation}{\mbox{\arabic{secnum}.\arabic{equation}}}}

We consider the  Boussinesq equations in the space time cylinder  $ \R^{2}\times (-1, 0)$
\begin{equation}
\begin{cases}
\partial _t v  + (v \cdot \nabla ) v  = e_1 \theta - \nabla p,\qquad  \nabla \cdot  u =0,
\\[0.3cm]
\partial _t \theta +  (v \cdot \nabla ) \theta =0,
\end{cases}
\label{bou}
\end{equation}
where $v=(v_1 (x,t), v_2 (x,t))$,  $(x,t)\in  \Bbb R^2\times (-\infty, 0)$.  This is an important equation modelling the dynamics of the heat convection in the atmospheric 
science(see e.g. \cite{maj1}). Moreover, it has essentially same structure as the axisymmetric 3D Euler equations off the axis(\cite{maj}). Therefore, the study of the system \eqref{bou} could provide us with information useful to understand the Euler equations. In \cite{cha1}  the first author of this paper proved the local well-posedness in standard Soblev space setting $H^m (\Bbb R^2)$, $m>2$, and also the following Beale-Kato-Majda\cite{bea} type (non)blow-up criterion is deduced; for any $m>2$
\begin{equation}
\lim\sup_{t\to 0^-} (\|v(t)\|_{H^m} +\|\theta(t)\|_{H^m} )<+\infty \quad \text{if and only if} \quad\int_{-1} ^0 \|\nabla \theta(t)\|_{L^\infty} dt<+\infty.
\label{thetaa}
\end{equation}
See also \cite{tan, dan} for the other forms of criterion,  using different functional setting, while a special type of scenario of singularity is excluded in \cite{cor}.
 We note the following scaling property of the system \eqref{bou}; it is invariant under the transform
$$(v(x,t), \theta (x,t))\mapsto \left(\lambda ^{\alpha} v(\lambda x, \lambda ^{\alpha+1} t), \lambda ^{2\alpha+1} \theta(\lambda x, \lambda ^{\alpha+1} t)\right)\quad \forall (x,t)\in \Bbb R^2 \times \Bbb R$$
for all $\lambda >1$, $\alpha >-1$. This leads us to the following natural definition.
\begin{defin}Let  $(v,\theta)\in C([-1, 0); W^{2,p} (\Bbb R^2))$, $p>2$,  be a local in time classical solution of \eqref{bou}, which blows up at $t=0$. We say it is of {\em Type I with respect to $v$},  if
$$
\sup_{-1<t<0} (-t) \|\nabla v (t)\|_{L^\infty} < +\infty,
$$
while we say it is of {\em Type I with respect to $\theta$}, if 
$$
\sup_{-1<t<0} (-t)^2 \|\nabla \theta (t)\|_{L^\infty} < +\infty.
$$
\end{defin}
 The aim of the present paper is to 
exclude a possible Type I blow-up at time $t=0$.    For the definition of  Besov space $\dot{B}^{0}_{\infty, \infty} (\Bbb R^n)$,  used in the theorem below,  see Section 2.

\begin{thm}
\label{thm1.1}
Let $ (v, \theta )$ be a solution of \eqref{bou} in $ \R^{2}\times (-1,0)$, which is regular in $ [-1,0)$. Furthermore, we assume that $ v(-1), \theta(-1) \in W^{2,\,p}(\R^{2})$ for $p>2$, and   at least one of the following conditions,
\begin{itemize}
\item[(i)]
$$ \lim\sup_{t\to 0^-} (-t)\|\nabla v(t)\|_{L^\infty} <2.
$$
\item[(ii)]
$$
 \intl_{-1}^{0} (-t) \| \nabla \theta(t)\|_{ L^\infty}   dt <+\infty. 
$$
\item[(iii)] 
$$
\int_{-1}^0 \|\o(t)\|_{\dot{B}^0_{\infty, \infty} }dt  +\int_{-1}^0(-t) \|\nabla \theta(t)\|_{\dot{B}^0_{\infty, \infty} }dt <+\infty.
$$
\item[(iv)]  
$$
\int_{-1}^0 \|\o(t)\|_{\dot{B}^0_{\infty, \infty} }dt+\sup_{-1<t<0} (-t)^2\|\nabla \theta(t)\|_{L^\infty } <+\infty.$$
 \end{itemize}
Then both $ v$ and $ \theta $ belong to $ L^\infty(-1, 0; W^{2,\, p}(\R^{2}))$. 

\end{thm}  
\begin{rem}  
(a)  In \cite{cha1} (see also \cite{hou}) it is proved that if  a solution to the 3D Euler equations on $\Bbb R^3\times  [-1, 0)$ satisfies $ v(-1) \in W^{2,\,p}(\R^{2})$ and
\begin{equation} 
\label{euler}\lim\sup_{t\to 0^-} (-t)\|\nabla v(t)\|_{L^\infty} <1.
\end{equation}
Then, $v\in L^\infty(-1, 0; W^{2,\, p}(\R^{2}))$. Note that the condition (i) is a relaxed version of  \eqref{euler}. 
 It is also interesting to notice that in a recent paper\cite{elg} Elgindi and Jeong  constructed 
explicitly blowing up solution,   which has linear growth at spatial infinity, and  is defined in a domain $D\subset \Bbb R^2$ with a corner. The  solution satisfies 
$$\lim_{t\to 0^-}\left\{(-t)\|\nabla v(t)\|_{L^\infty (D)} +(-t)^2 \|\nabla \theta (t)\|_{L^\infty (D)}\right\} <+\infty, $$
and the blow-up happens at every point in $D$.
\\
\ \\
(b) The main novelty in  the conditions (ii) and (iii) is the extra factor $(-t)$ in the integral of the norms of $|\nabla \theta(t)|$. This factor makes  the integral $ \int_{-1} ^0(-t) \|\nabla \theta(t)\|_{L^\infty} dt$ scaling invariant quantity, while  the stronger  integral $  \int_{-1} ^0\|\nabla \theta(t)\|_{L^\infty} dt$  is not. Similar remark holds for  $\int_{-1}^0 \|\nabla \theta(t)\|_Xdt $ with $X=BMO$ or $\dot{B}^0_{\infty, \infty}$.\\
\ \\
(c) As far as the authors know it is still an open question if the regularity of the system \eqref{bou} is guaranteed only by the vorticity integral condition, say
\begin{equation}  
\int_{-1}^0 \|\omega (t)\|_{L^\infty} dt <+\infty.
\label{vort}
\end{equation}  The above theorem with the condition (iv) says that  if \eqref{vort}  holds, then any singularity, which is of Type I  with respect to $\theta$  is excluded.
\end{rem}

At this moment we could not omit the vorticity integral in the condition (iv) above, but if we modify  Type I condition on $\nabla \theta$ logaritmmically as well as imposing the smallness, then this is possible as follows.

\begin{thm}
\label{thm1.2}
Let $ (v, \theta )$ be a solution of \eqref{bou} in $ \R^{2}\times (-1,0)$ which is regular in $ [-1,0)$. Furthermore, we assume that $ v(-1), \theta(-1) \in W^{2,\, p_0}(\R^{2})$ for some $ p_0 >2$. 
 There  exists $ \var >0$ depending only on $ p_0$, such that if 
\begin{equation}
\limsup_{t \to 0^-} (-t)^2  \log (-1/t) \| \nabla \theta (t)\|_{ L^\infty} \le \var,   
\label{type1e}
\end{equation}
then both $ v$ and $ \theta $ belong to $ L^\infty(-1, 0; W^{2,\, p_0}(\R^{2}))$. 
\end{thm}  

\section{Proof of the Main Theorems}
\label{sec:-5}
\setcounter{secnum}{\value{section} \setcounter{equation}{0}
\renewcommand{\theequation}{\mbox{\arabic{secnum}.\arabic{equation}}}}

We introduce the space $\dot{B}^{0}_{\infty, \infty} (\Bbb R^n)$ below.
Let $\varphi \in \mathcal{S}$, where  $\mathcal{S}$ is  the the Schwartz class of rapidly decreasing functions, and let $\hat{\varphi}$ be its Fourier transform, defined by
$\hat{\varphi}(\xi)= \int_{\Bbb R^n} e^{-2\pi i x\cdot \xi } \varphi(x)dx.$ Then, we consider $\varphi$ satisfying the following conditions
$$Supp \,\hat{\varphi}\subset \{ \xi \in \Bbb R^n\, |\, \frac12\leq |\xi|\leq 2\},\quad
\hat{\varphi}\geq c>0 \,\quad \text{if} \quad\frac23 <|\xi|<\frac32, \quad \text{and}\quad\sum_{j\in \Bbb Z}  \hat{\varphi}_j (\xi)=1, 
$$
where we defined $ \hat{\varphi}_j= \hat{\varphi}(2^{-j}\xi )$. Construction  of  the sequence $\{ \varphi_j\}_{j\in \Bbb Z}$ is well-known(see e.g. \cite{che}).  Then, we say  $f\in \dot{B}^{0}_{\infty, \infty} (\Bbb R^n)$  if and only if
$ \sup_{j\in \Bbb R^n} \|\varphi_j \ast f \|_{L^\infty}:= \|f \|_{\dot{B}^{0}_{\infty, \infty} }<+\infty$.
The basic properties of $ \dot{B}^{0}_{\infty, \infty} (\Bbb R^n)$ useful for us are the followings.
\begin{itemize}
\item[(i)] Embedding properties:
\bb\label{emb}
L^\infty (\Bbb R^n)\hookrightarrow BMO(\Bbb R^n)\hookrightarrow\dot{B}^{0}_{\infty, \infty} (\Bbb R^n),
\ee
\item[(ii)] The logarithmic Sobolev inequality,
\bb
\|f\|_{L^\infty}\leq c (1+\|f\|_{\dot{B}^{0}_{\infty, \infty} } \log(e+ \|f\|_{W^{s,p}} )), \quad s>n/p.
\ee
where the constant $c=:c_{ls}$ depends on $s$ and $p$.\\
\item[(iii)] Boundedness on the Calderon-Zygmund operators, in particular applying to the Bio-Savart formular one has
\bb \|\nabla v\|_{\dot{B}^{0}_{\infty, \infty} }\leq c \|\omega\|_{\dot{B}^{0}_{\infty, \infty} },
\ee
where $(v, \omega)$ satisfies $\nabla \cdot v=0, \nabla \times v=\omega$.
\end{itemize}
 \ \\
\noindent{\bf Proof of Theorem\ref{thm1.1} :} \\
\noindent\underline{{\it Proof for (iii): }} 
Let $q>2$. We  apply the operator $ \partial _i $ to the vorticity equation, multiplying  the resultant equation by 
$ \partial _i \omega | \nabla \omega |^{q-1} $, and integrating it over $\Bbb R^2$.  Then, after the integration by part  and  using the H\"{o}lder inequality, we are led to
\begin{align}
&\frac{d}{dt}\| \nabla \omega\|_{L^q} \le  \| \nabla v\|_{L^\infty} \| \nabla \omega \|_{L^q} + \| \nabla ^2 \theta \|_{L^q}\cr
&\qquad=  \| \nabla v\|_{L^\infty} \| \nabla \omega \|_{L^q} + (-t)^{-1} (-t)\| \nabla ^2 \theta \|_{L^q}. 
\label{dieq1}
\end{align} 

\hspace{0.5cm}
Next, we apply the operator $ \partial _i \partial _j$ to both sides of the   $ \theta $ equation,  multiply both sides the by  $ \partial _i\partial_j \theta  | \nabla ^2 \theta |^{q-2}$, and sum over $i,j=1,2,3$, and the integrate it over $\Bbb R^2$. 
This, applying  the integration by part and the H\"{o}lder inequality, yields the following  inequality
\begin{equation}
\frac{d}{dt}\| \nabla^2 \theta  \|_{L^q}  \le  
 2 \| \nabla v\|_{L^\infty}\| \nabla ^2 \theta\|_{L^q}  + \| \nabla \theta\|_{L^\infty}\| \nabla ^2 v\|_{L^q}.    
 \label{dieq3}
\end{equation}

Multiplying both sides of \eqref{dieq3}  by $ (-t)$,  we see that 
\begin{align}
&\frac{d}{dt}(-t)\| \nabla^2 \theta  \|_{L^q} +\| \nabla^2 \theta  \|_{L^q} \cr
&\qquad \le 2 \| \nabla v\|_{L^\infty}(-t)\| \nabla ^2 \theta\|_{L^q}  + (-t)\| \nabla \theta\|_{L^\infty}\| \nabla ^2 v\|_{L^q}\cr
 &\qquad\leq   2  \| \nabla v\|_{L^\infty}(-t)\| \nabla ^2 \theta\|_{L^q}  + c_{cz}(-t)\| \nabla \theta\|_{L^\infty}\| \nabla \omega\|_{L^q}\label{dieq4}
\end{align}

Now define 
\[
\Psi (t) := \|  \nabla \omega\|_{ L^q}  + (-t) \| \nabla ^2 \theta \|_{ L^q},\quad  t\in (-1,0).
\]
Adding the last two inequalities \eqref{dieq1} and \eqref{dieq4}, we are led to 
\begin{align}
&\Psi' \le \Big(2\| \nabla v(t)\|_{ L ^\infty}  +  (-t)^{-1}+  c_{cz}(-t)\| \nabla \theta (t)\|_{ L^\infty} \Big)\Psi .
\label{5.4}
\end{align}

\hspace{0.5cm}
By means of the logarithmic Sobolev embedding,   we find 
\begin{align}
\| \nabla v(t)\|_{ L^\infty}
& \le c\left\{1+\| \nabla v(t)\|_{\dot{ B}^0_{\infty, \infty}} \log (e+ \| \nabla ^2 v(t)\| _{L^ q}) \right\}
\cr
&\le  c\left\{1+ \| \omega(t)\|_{\dot{ B}^0_{\infty, \infty}}  \log (e+ \| \nabla  \omega(t) \| _{L^ q})  \right\}
\cr
&\le  c \left\{1+\| \omega(t)\|_{\dot{ B}^0_{\infty, \infty}}  \log (e+ \Psi (t)) \right\}. 
\label{5.6}
\end{align}
Similarly,
\bb\label{theta}
 \|\nabla \theta\|_{L^\infty}  \le  c\left\{1+ \| \nabla \theta(t)\|_{\dot{ B}^0_{\infty, \infty}}  \log (e+ \Psi (t)) \right\}. 
 \ee

Inserting \eqref{5.6} and \eqref{theta} into \eqref{5.4},  it follows 
 \begin{align}
\Psi' &\le \left\{c \left[1+  (\| \omega(t)\|_{\dot{ B}^0_{\infty, \infty}} +  (-t)\| \nabla \theta(t)\|_{\dot{ B}^0_{\infty, \infty}} )  \log (e+ \Psi (t))\right]+ (-t)^{ -1}\right\}\Psi(t).
\label{5.7}
\end{align}

Setting $ y(t)= \log (e+ \Psi (t))$, we infer  from \eqref{5.7} the differential inequality  
\begin{equation}
y' \le  c a(t)y + c(-t)^{ -1} ,\qquad  a(t)= \| \omega(t)\|_{\dot{ B}^0_{\infty, \infty}} +  (-t)\| \nabla \theta(t)\|_{\dot{ B}^0_{\infty, \infty}} \label{5.8}
\end{equation}
which can be solved as
\begin{align}
\label{5.9}
y(t)=&\log (e+ \Psi (t))\cr
&\le y(t_0) e^{ c\int_{t_0}^{t} a(s) ds} +  c \intl_{t_0}^{t}  (-s)^{ -1} e^{ c\int^{t }_{s} a(\tau ) d\tau }  ds 
\end{align}
 We now choose $t_0$ so that 
$e^{c\int_{t_0} ^0 a(s)ds} <2$. Then, \eqref{5.9} implies
\bb\label{log} \log (e+ \Psi (t)) \leq c \log (e+ \Psi (t_0))+c\log (-1/t) \qquad \forall t\in (t_0, 0),
\ee
where $c>2$ is another constant.
From $\theta$-equation we have immediately
\bb\label{theq}
\frac{\partial}{\partial t} |\nabla \theta| +(v\cdot\nabla )|\nabla \theta|\leq	  |\nabla v||\nabla \theta |.
\ee
Let $ t\in (-1,0)$ be arbitrarily chosen but fixed. Let $ x_0\in \R^{2}$. By $ X(x_0,t)$ we denote the trajectory of the particle 
which is located at $ x_0$ at time $ t=t_0 $, defined by the following ODE
\begin{equation}
\frac{d X(x_0, t)}{d  t} = v(X(x_0 ,t), t)\quad  \text{ in }\quad  [-1,0),\quad  X(x_0 , t_0) = x_0.  
\label{5.1}
\end{equation}
The Lipschitz continuity of   $ v(s)$ in $ \R^{2}$ for all $ s\in (-1,0)$ ensures the existence and uniqueness a  solution to  \eqref{5.1} in $[-1,0)$. 
Then, \eqref{theq} can be written as
\bb\label{theq0}
\frac{\partial}{\partial t} |\nabla \theta(X(x_0, t),t)| \leq	  |\nabla v(X(x_0, t),t)||\nabla \theta(X(x_0, t),t) |,
\ee
which can be integrated along the trajectories as
$$
 |\nabla \theta(X(x_0, t),t)| \leq  |\nabla \theta(x_0)|  \exp\left( \int_{t_0} ^t |\nabla v( X(x_0, s), s)|ds\right).
 $$
 Therefore, we estimate, using \eqref{log} as
 \begin{align}\label{theqq}
\|\nabla \theta (t)\|_{L^\infty} & \leq \|\nabla \theta(t_0)\|_{L^\infty} \exp\left(\int^t_{t_0} \|\nabla v\|_{L^\infty} ds\right)\cr
  \leq &\|\nabla \theta(t_0)\|_{L^\infty} \exp\left(c\int^t_{t_0} \left\{ \|\omega(s)\|_{\dot{ B}^0_{\infty, \infty}} \left[
  \log (e+ \Psi (t_0))+\log (-1/s)\right] +1\right\}ds\right)\cr
   \leq &\|\nabla \theta(t_0)\|_{L^\infty} \exp\left(c\left\{
  \log (e+ \Psi (t_0))+\log (-1/t)\right\}\int^t_{t_0} \|\omega(s)\|_{\dot{ B}^0_{\infty, \infty}} ds +c(t-t_0)\right) \end{align}
Choosing $t_0\in (-1,0)$ so  that
$$
c\int^0_{t_0} \|\omega(s)\|_{\dot{ B}^0_{\infty, \infty}} ds <\frac{1}{2},
$$
we deduce from  \eqref{theqq} that
$$
\|\nabla \theta (t)\|_{L^\infty}\leq \|\nabla \theta(t_0)\|_{L^\infty}  (e+ \Psi(t_0))^{c} e^{c}(-t)^{-\frac12} \quad \forall t\in (t_0, 0).
$$
Therefore, $\int_{-1} ^0 \|\nabla\theta \|_{L^\infty} dt <+\infty.$ Applying the well-known blow-up criterion in \cite{cha2},
we obtain the desired result.\\
\ \\
\noindent{\it \underline{Proof  for (iv) :} }  Under the hypothesis of (iv) \eqref{5.4} is replaced by
\begin{align}
&\Psi' \le \Big(2\| \nabla v(t)\|_{ L ^\infty}  +  c(-t)^{-1}\Big)\Psi ,
\end{align}
and the remaing part of the proof is the same as in (iii).\\
\ \\
\noindent{\it \underline{Proof  for (ii):} }
Applying curl to the velocity equation in \eqref{bou},   we obtain 
\begin{equation}
\partial _t \omega + v\cdot \nabla \omega = - \partial _2 \theta\quad  \text{ in}\quad  \R^{2}\times [-1, 0),  
\label{vorteq}
\end{equation}
where $ \omega = \partial _1 v_2 - \partial _2 v_1$.

\hspace{0.5cm}
Using the particle trajectories(with $X(x_0, -1)= x_0$) as the above,  we have from \eqref{vorteq} 
\begin{equation}
\frac{d}{dt} |\omega (X(x_0, t),t)| \leq |\partial _2\theta (X(x_0, t),  t)|\quad \text{ in}\quad  [-1,0),
\label{5.1a}
\end{equation}
which implies that
\bb\label{5.1aa}
\|\omega(s)\|_{L^\infty} \leq \|\omega(-1)\|_{L^\infty} +\int_{-1} ^s \|\partial_2 \theta(\tau)\|_{L^\infty} d\tau.
\ee
Integrating both sides of \eqref{5.1aa} over $[-1, t)$, $t\in (-1, 0)$ with respect to $s$,  and applying integration by parts, we get
\begin{align*}
\int_{-1} ^t \|\omega(s)\|_{L^\infty} ds &\le (1+t) \|\omega(-1)\|_{L^\infty} + \int_{-1} ^t\int_{-1} ^s \|\partial_2 \theta(\tau)\|_{L^\infty} d\tau ds\\
&=(1+t) \|\omega(-1)\|_{L^\infty}  + \int_{-1} ^t\left\{ \frac{d}{ds} (s)\int_{-1} ^s \|\partial_2 \theta(\tau)\|_{L^\infty} d\tau \right\}ds\\
&= (1+t) \|\omega(-1)\|_{L^\infty} + \int_{-1} ^t (-s) \|\partial_2 \theta(s)\|_{L^\infty} ds +t \int_{-1} ^t \|\partial_2 \theta(s)\|_{L^\infty} ds\\
&\le  \|\omega(-1)\|_{L^\infty} + \int_{-1} ^t (-s) \|\partial_2 \theta(s)\|_{L^\infty} ds.
\end{align*}
Therefore, 
 \begin{align}
&\intl_{-1}^{t}  \| \omega (s)\|_{ L^\infty}  ds +   \intl_{-1}^{t}  (-s) \| \nabla \theta (s)\|_{ L^\infty}  ds\cr
&\qquad\qquad\leq\|  \omega (-1) \|_{ L^\infty}
 + 2 \intl_{-1}^{0} (-s)\| \nabla \theta (s)\|_{ L^\infty}   ds<+\infty. 
\label{bkm}
\end{align}
Therefore, from the embedding \eqref{emb}  the condition (iii) is satisfied.\\
\ \\
\noindent{\it \underline{Proof  for (i) :} }
By hypothesis (i)  there exists $t_0\in (-1, 0)$ and $\delta >0$ such that
$$ \sup_{t_0<t<0} (-t )\|\nabla v(t)\|_{L^\infty} \leq 2-\delta.$$
Multiplying \eqref{theq} by $-\tau$, we have
$$
\frac{\partial}{\partial  \tau}  ( (-\tau) |\nabla \theta| ) +|\nabla \theta| +(v\cdot \nabla) (-\tau |\nabla \theta|) 
 \le (-\tau) |\nabla v||\nabla \theta| \le (2-\delta ) |\nabla \theta|,
$$
which after integration over $(t_0, s)$ along the trajectory gives
$$
(-s) |\nabla \theta (X(x_0, s),s)|\le (-t_0) |\nabla \theta(x_0,t_0)|+(1-\delta ) \int_{t_0} ^ s \|\nabla \theta (\tau) \|_{L^\infty} d\tau.
$$
Let $t\in (t_0, 0)$. Then, for all $ s\in (t_0, t)$ we have
$$
(-s) \|\nabla \theta (s)\|_{L^\infty} \le (-t_0) \|\nabla \theta (t_0)\|_{L^\infty} +(1-\delta) \int_{t_0} ^s \|\nabla \theta (\tau) \|_{L^\infty}d\tau.
$$
Integrating the both sides of the above over $(t_0, t)$, and integrating by part, we get
\begin{align*}
&\int_{t_0} ^t (-s) \|\nabla \theta (s)\|_{L^\infty} ds \le (-t_0) (t-t_0) \|\nabla \theta (t_0)\|_{L^\infty} +(1-\delta) \int_{t_0} ^s \|\nabla \theta (\tau) \|_{L^\infty} d\tau\\
 & \le (-t_0)^2 \|\nabla \theta (t_0)\|_{L^\infty} +(1-\delta)\left\{ \left(s\|\nabla \theta (s) \|_{L^\infty} -t_0\|\nabla \theta (t_0) \|_{L^\infty}\right)  - \int_{t_0} ^s \tau \|\nabla \theta (\tau) \|_{L^\infty} d\tau \right\}\\
 &\le   (-t_0)^2 \|\nabla \theta (t_0)\|_{L^\infty}  +(1-\delta) ( -t_0) \|\nabla \theta (t_0) \|_{L^\infty}  +(1-\delta)\int_{t_0} ^t (- \tau) \|\nabla \theta (\tau) \|_{L^\infty} d\tau.
   \end{align*}
which implies
$$
\delta \int_{t_0} ^t (-s) \|\nabla \theta (s)\|_{L^\infty} ds \le (-t_0)^2 \|\nabla \theta (t_0)\|_{L^\infty}  +(1-\delta) ( -t_0) \|\nabla \theta (t_0) \|_{L^\infty},
$$
Passing $t\to 0^-$, we obtain finally
$$
\delta \int_{t_0} ^0 (-s) \|\nabla \theta (s)\|_{L^\infty} ds \le (-t_0)^2 \|\nabla \theta (t_0)\|_{L^\infty}  +(1-\delta) ( -t_0) \|\nabla \theta (t_0) \|_{L^\infty} <+\infty 
$$
and  the condition (ii) is satisfied.  \hfill \Beweisende   \\
\ \\
\ \\
\noindent{\bf Proof of Theorem\,\ref{thm1.2} :} From  \eqref{type1e}, we find $ t_0\in (-e^{ -2},0) $ such that 
\begin{equation}
\| \nabla \theta (s)\|_{ L^\infty}  \le  \frac{\var }{(-s)^2 \log(-1/s)} \quad  \forall\,s\in [t_0, 0).
\label{6.0}
\end{equation}
The inequality \eqref{5.1a}, following the argument of the proof for (ii), and combined with  \eqref{6.0} yields
\begin{align}
\| \omega (t)\|_{ L^\infty} &\le \| \omega (t_0)\|_{ L^\infty} +  \intl_{t_0}^{t} \| \nabla \theta (s)\|_{ L^\infty}   ds
\cr
&\le \| \omega (t_0)\|_{ L^\infty} + \var  \intl_{t_0}^{t} \frac{1}{(-s)^2 \log(-1/s)}   ds
\cr
&\le \| \omega (t_0)\|_{ L^\infty} + 2\var  \intl_{t_0}^{t} \frac{\log(-1/s)-1}{((-s) \log(-1/s))^2}   ds
\cr
&= \| \omega (t_0)\|_{ L^\infty} + \frac{2\var }{(-t)\log(-1/t)}-\frac{2\var }{(-t_0)\log(-1/t_0)}
\cr
& \le \| \omega (t_0)\|_{ L^\infty} + \frac{2\var }{(-t)\log(-1/t)},
\label{6.1}
\end{align}
where we used the fact that $\log (-1/s) \leq 2\log (-1/s) -2$ for all $s\in (-e^{-2}, 0)$ in the third inequality.
We now define  $ \var >0$ as follows 
\begin{equation}
\var := \frac{1}{4 \max\{ c_{ ls}, c_{ cz}\}}. 
\label{condeps}
\end{equation}
Then from \eqref{5.4}  combined with \eqref{5.6} together with \eqref{6.0} and \eqref{6.1}  we find
 \begin{align}
y' &\le c_{ ls} \| \omega\|_{\infty} y  +    (-t)^{ -1}+  c_{ cz}(-t) \| \nabla \theta(t) \|_{ L^\infty}
\cr
 &\le \Big( \frac{2\vare c_{ls}}{(-t)\log(-1/t)} + c_{ ls}  \| \omega (t_0)\|_{ L^\infty} \Big) y +  (-t)^{ -1} +\vare c_{cz} (-t)^{-1} \log(-1/t)\cr
 &  \le \Big( \frac{1}{2(-t)\log(-1/t)} + c_{ ls}  \| \omega (t_0)\|_{ L^\infty} \Big) y + \frac{5}{4} (-t)^{ -1}     \quad  \text{ in}\quad  (t_0,0),
\label{6.2}
\end{align}
where 
\[
y(t) = \log (e+ \Psi (t)),\quad  \Psi (t) := \|  \nabla \omega(t)\|_{ L^{ p_0}}  + (-t)
\| \nabla ^2 \theta(t) \|_{ L^{ p_0}},\quad  t\in (-1,0).
\]
Integrating \eqref{6.2}, we obtain 
\begin{align}
y(t) \le y(t_0) e^{  \intl_{t_0}^{t} a(s)  ds} +  \frac{5}{4}\intl_{t_0}^{t} (-s)^{ -1}    e^{  \intl_{s}^{t} a(\tau )  d\tau }ds,
\label{6.4}
\end{align}
where we set
\begin{align*}
a(t) = \frac{1}{2(-t)\log(-1/t)} + c_{ ls}  \| \omega (t_0)\|_{ L^\infty} . 
\end{align*}
Applying integration by parts, we infer 
\begin{align*}
&\intl_{t_0}^{t} (-s)^{ -1}    e^{  \intl_{s}^{t} a(\tau )  d\tau }ds 
\\
&\quad = \intl_{t_0}^{t} \frac{d}{ds} \log(-1/s)  e^{  \intl_{s}^{t} a(\tau )  d\tau }ds
\\
&\quad  = \log(-1/t) - \log(-1/t_0)   e^{  \intl_{t_0}^{t} a(\tau )  d\tau } + 
\intl_{t_0}^{t} \log(-1/s)  a(s) e^{  \intl_{s}^{t} a(\tau )  d\tau }ds
\\
& \quad \le \log(-1/t) + \frac{1}{2}  \intl_{t_0}^{t} (-s)^{ -1}    e^{  \intl_{s}^{t} a(\tau )  d\tau }ds  
+ c_{ ls} \| \omega (t_0)\|_{ L^\infty}\intl_{t_0}^{t} \log(-1/s)    e^{  \intl_{s}^{t} a(\tau )  d\tau }ds.
\end{align*}
Absorbing the second term on the right hand side into the left, one has
\begin{align*}
&\intl_{t_0}^{t} (-s)^{ -1}    e^{  \intl_{s}^{t} a(\tau )  d\tau }ds \cr
\quad&\le 2\log(-1/t) +2c_{ ls} \| \omega (t_0)\|_{ L^\infty}\intl_{t_0}^{t} \log(-1/s)    e^{  \intl_{s}^{t} a(\tau )  d\tau }ds.
\end{align*}
Calculating 
\begin{align}
e^{  \intl_{s}^{t} a(\tau )  d\tau } &=  e ^{c_{ ls} \| \omega (t_0)\|_{ L^\infty} }  
e ^{ \frac{1}{2}\intl_{s}^{t} \frac{1}{(-\tau) \log(-1/\tau ) }  d\tau} 
\leq e ^{c_{ ls}  \| \omega (t_0)\|_{ L^\infty} }  \left\{\log (-1/t )\right\}^{ \frac{1}{2}}
\label{6.5}
\end{align}
for all $ s\in [t_0, 0)$, we obtain from the above inequality 
\begin{align*}
&\intl_{t_0}^{t} (-s)^{ -1}    e^{  \intl_{s}^{t} a(\tau )  d\tau }ds  
\\
&\quad  \le 2\log(-1/t) +2c_{ls}\| \omega (t_0)\|_{ L^\infty} e ^{c_{ ls}  \| \omega (t_0)\|_{ L^\infty} }  \{\log (-1/t )\}^{\frac12} \intl^{0}_{-1}\log(-1/s)  ds 
  \cr
 &\quad  \le 2\log(-1/t)+ c \{ \log(-1/t)\}^{ \frac{1}{2}},
\end{align*}
where $ c=\const$ is independent on $ t$. Estimating the second term in \eqref{6.4} by the estimate we have just obtained and the first term by  \eqref{6.5} for $ s=t_0$, 
we arrive at  
\[
y(t) \le \frac{5}{2} \log(-1/t)+ c \{ \log(-1/t)\}^{ \frac{1}{2}}  \quad  \forall\,t\in [t_0, 0),
\]
for some constant independent of $ t$.  Accordingly, there exists $ t_1\in (t_0, 0)$ such that 
\begin{equation}
y(t) \le 3 \log(-1/t) \quad  \forall\,t\in [t_1, 0),
\label{6.6}
\end{equation}

\hspace{0.5cm}
By the aid of the logarithmic Sobolev embedding inequality, and  observing \eqref{6.1} together with \eqref{6.6} 
and \eqref{condeps}, we see that 
 for all $ t\in [t_1, 0)$
 \begin{align}
 \| \nabla v(t)\|_{ L^\infty} &\le  c_{ ls} \| \omega(t) \|_{ L^\infty} y(t) +c_{ls}\le 3c_{ls} \| \omega (t_0) \|_{ L^\infty} \log(-1/t)+  6\var  c_{ ls}    (-t)^{ -1}+c_{ls} \cr
 &\le 3c_{ ls} 
 \| \omega (t_0) \|_{ L^\infty}  \log(-1/t)+  \frac{3}{2}  (-t)^{ -1} +c_{ls} .\label{6.12}
 \end{align}
Thus,  
\begin{equation}
\limsup_{ t\to 0^-} (-t) \| \nabla v(t)\|_{ L^\infty} \le \frac{3}{2} < 2. 
\label{6.13}
\end{equation}  
Applying  Theorem \ref{thm1.1} (i), we get the assertion of the theorem.  \hfill \Beweisende \\
\ \\
$$\mbox{\bf Acknowledgements}$$
Chae was partially supported by NRF grants 2016R1A2B3011647, while Wolf has been supported 
supported by NRF grants 2017R1E1A1A01074536. The authors thank to the anonymous referees for many valuable suggestions. The authors declare that they have no conflict of interest. 
  \end{document}